\documentclass{article}
\usepackage{fullpage}
\usepackage[final]{hyperref}
\usepackage{nameref}
\usepackage[round]{natbib}
\usepackage{tabularx}
\usepackage{booktabs}
\usepackage{setspace}
\usepackage{multirow}
\usepackage{makecell}
\usepackage{graphicx}
\usepackage{delarray}
\usepackage{xcolor}
\usepackage{colortbl}
\usepackage{amssymb}
\usepackage{amsmath}
\usepackage{mathrsfs}
\usepackage{mathtools}
\usepackage{nicefrac}       
\usepackage{tikz}
\usepackage{xspace}
\xspaceaddexceptions{[]\{\}}
\usepackage{algorithm}
\usepackage{algorithmic}
\usepackage{mdframed}
\usepackage{bm}
\usepackage{tablefootnote}

\newcommand{\cD}{{\mathcal D}}

\newtheorem{assumption}{Assumption}
\newtheorem{theorem}{Theorem}
\newtheorem{lemma}{Lemma}

{\hspace*{\fill}$\Box$\par}
{\hspace*{\fill}$\Box$\par\vspace{4mm}}
{\hspace*{\fill}$\Box$\par}

\newcommand{\framed}[1]{
	\begin{mdframed}[linecolor=black,backgroundcolor=gray!10]
		#1
	\end{mdframed}
}

\hypersetup{
	colorlinks=true,       
	linkcolor=blue,        
	citecolor=blue,        
	filecolor=magenta,     
	urlcolor=blue
}

\newcommand*{\qedb}{\hfill\ensuremath{\square}}   

\newcommand{\head}[1]{\noindent{{\bf #1:}}}
\newcommand{\topic}[1]{\vspace{2mm}\noindent{{\bf #1:}}}

\definecolor{bgcolor}{rgb}{0.66,0.88,1.00}
\definecolor{bgcolor2}{rgb}{0.66,0.88,0.50}

\definecolor{harvardcrimson}{rgb}{0.79, 0.0, 0.09}
\definecolor{darkviolet}{rgb}{0.58, 0.0, 0.83}
\newcommand{\cfinite}[1]{\textcolor{harvardcrimson}{#1}}
\newcommand{\conline}[1]{\textcolor{darkviolet}{#1}}

\newcommand{\E}{{\mathbb{E}}}
\newcommand{\R}{{\mathbb{R}}}

\newcommand{\ns}[1]{\| #1 \|^2}

\newcommand{\n}[1]{\| #1 \|}

\newcommand{\hx}{\widehat{x}}

\newcommand{\page}{{\sf PAGE}\xspace}

\newcommand{\fgap}{{\Delta_0}}

\newcommand{\xt}{x^t}
\newcommand{\xtn}{x^{t+1}}

\newcommand{\gt}{g^t}
\newcommand{\gtn}{g^{t+1}}

\newcommand{\EB}[1]{\E\left[ #1 \right]}
\newcommand{\newl}{\notag \\ &\qquad }

\begin{document}

\title{\bf \Large A Short Note of PAGE:\\
	Optimal Convergence Rates for Nonconvex Optimization}
\author{Zhize Li
			\\KAUST}

\date{\today}
\maketitle

\begin{abstract}
	In this note, we first recall the nonconvex problem setting and introduce the optimal \page algorithm \citep{li2021page}. 
	Then we provide a simple and clean convergence analysis of \page for achieving optimal convergence rates.
	Moreover, \page and its analysis can be easily adopted and generalized to other works.
	We hope that this note provides the insights and is helpful for future works.
\end{abstract}


\section{Problem Setting}
\label{sec:prob}

We consider the nonconvex optimization problem
$
\min_{x\in \R^d}   f(x).
$
The nonconvex function $f$ has the following two forms:
\begin{enumerate}
	\item \cfinite{\emph{Finite-sum}} form
		\begin{align}
		\cfinite{\min_{x\in \R^d} \bigg\{f(x) := \frac{1}{n}\sum_{i=1}^n{f_i(x)} \bigg\}}, \label{prob:finite}
		\end{align}
		where functions $f_i$s are differentiable and possibly nonconvex, e.g., there are $n$ data samples and $f_i$ is a nonconvex loss on data $i$;

	\item \conline{\emph{Online}} form
		\begin{align}
		\conline{\min_{x\in \R^d}  \bigg\{ f(x) := \E_{\zeta\sim \cD}[F(x,\zeta)] \bigg\}}, \label{prob:online}
		\end{align}
		where $F(x,\zeta)$ is also differentiable and possibly nonconvex, e.g., in the online/streaming case, data is drawn from an unknown distribution $\cD$.
\end{enumerate}
For notational convenience, we will simply denote the online form as the finite-sum form via letting $f_i(x) := F(x, \zeta_i)$ and treating $n$ as a very large value or even infinite.

\vspace{2mm}
Now we define the following standard assumptions. 

\begin{assumption}[Average $L$-smoothness]\label{asp:smooth}
	The function $f$ is average $L$-smooth if $\exists L >0$, 
	\begin{equation}\label{eq:avgsmooth}
	\E_i[\ns{\nabla f_i(x) - \nabla f_i(y)}]\leq L^2 \ns{x-y}, ~ \forall x,y \in \R^d.
	\end{equation}
\end{assumption}

\begin{assumption}[Bounded variance]\label{asp:bv}
	The stochastic gradient has bounded variance if $\exists \sigma >0$,
	\begin{equation}\label{eq:bv}
	\E_{i}[\ns{\nabla f_i(x)-\nabla f(x)}]\leq \sigma^2, \quad \forall x \in \R^d.
	\end{equation}
\end{assumption}

Note that the average $L$-smoothness Assumption \ref{asp:smooth} implies $f$ is $L$-smooth (see Lemma 1 of \citealp{li2021page}). 
The Assumption \ref{asp:bv} usually is only needed for the online case \eqref{prob:online} since the full gradient  ($\nabla f(x)=\frac{1}{n}\sum_{i=1}^n{\nabla f_i(x)}$) may not be available (e.g., if $n$ is infinite).

\section{PAGE Algorithm}
\label{sec:algo}

In each iteration, the gradient estimator $g^{t+1}$ of \page is defined in  Line \ref{line:grad} of Algorithm \ref{alg:page}, 
which indicates that \page uses the vanilla minibatch SGD update with probability $p_t$, and reuses the previous gradient $g^{t}$ with a small adjustment (which  lowers the computational cost if $b' \ll b$) with probability $1-p_t$. 
In particular, the  $p_t \equiv 1$ case reduces to minibatch SGD, and to GD if we further set the minibatch size $b=n$.

\vspace{-1mm}
\begin{algorithm}[htb]
	\caption{ProbAbilistic Gradient Estimator (\page) \citep{li2021page}}
	\label{alg:page}
	\begin{algorithmic}[1]
		\REQUIRE ~
		initial point $x^0$, stepsize $\eta$, minibatch size $b$ and $b'$, probability $\{p_t\} \in (0,1]$
		\STATE $g^{0} = \frac{1}{b} \sum_{i\in I} \nabla f_i(x^0)$  \quad {\footnotesize ~//~$I$ denotes random minibatch samples with $|I|=b$} \label{line:sgd}
		\FOR {$t=0,1,2,\ldots$}
		\STATE $x^{t+1} = x^t - \eta g^t$  \label{line:update}
		\STATE $g^{t+1} = \begin{cases}
		\frac{1}{b} \sum \limits_{i\in I} \nabla f_i(x^{t+1}) &\text{with probability } p_t\\
		g^{t}+\frac{1}{b'} \sum \limits_{i\in I'} (\nabla f_i(x^{t+1})- \nabla f_i(x^{t})) &\text{with probability } 1-p_t
		\end{cases}$   \label{line:grad}
		\ENDFOR
		\ENSURE $\hx_T$ chosen uniformly from $\{x^t\}_{t\in [T]}$  
	\end{algorithmic}
\end{algorithm}

\vspace{-5mm}
\subsection{Optimal Convergence rates of PAGE}

In Theorem \ref{thm:main}, we show that a parameter choice of \page can lead to optimal convergence rates, matching the lower bound \cfinite{$\Omega\big(n+\frac{\sqrt{n}}{\epsilon^2}\big)$} (Theorem 2 of \citealp{li2021page}) and \conline{$\Omega\big(b+\frac{\sqrt{b}}{\epsilon^2}\big)$} (Corollary 5 of \citealp{li2021page}).
\framed{
\begin{theorem}\label{thm:main}
	Let stepsize $\eta \leq \frac{1}{L\left(1+\sqrt{\frac{1-p}{pb'}}\right)}$ and probability $p_t\equiv \frac{b'}{b+b'}$.
	Then \page (Algorithm \ref{alg:page}) can find an $\epsilon$-approximate solution, i.e., 
	$${\E[\n{\nabla f(\hx_T)}]\leq \epsilon}.$$ 
	We distinguish the following two cases:
	\begin{enumerate}
		\item (\cfinite{Finite-sum} case) Under Assumption \ref{asp:smooth}, let minibatch size $b=n$ and any $b' \leq \sqrt{b}$, then the number of iterations can be bounded by 
		\vspace{-1mm}
		$$ T= \frac{2 L \fgap }{\epsilon^2} \left( 1+ \sqrt{\frac{1-p}{pb'}}\right) \leq \frac{4L\fgap \sqrt{n}}{\epsilon^2 b'},
		$$
		and the number of stochastic gradient computations (i.e., gradient complexity) is 
		 $$ \#\mathrm{grad} =b+ T\left(pb+(1-p)b'\right) \leq  n + \frac{8L\fgap \sqrt{n}}{\epsilon^2} = \cfinite{O\bigg(n+\frac{\sqrt{n}}{\epsilon^2}\bigg)},
		 $$
		where $\fgap:=f(x^0)-\min_{x\in \R^d}f(x)$, and the first $b$ in $\#\mathrm{grad}$ is due to $g^0$ (Line \ref{line:sgd} in Algorithm \ref{alg:page}).

		\item (\conline{Online} case) Under Assumptions \ref{asp:smooth} and \ref{asp:bv}, let minibatch size $\conline{b=\min\big\{\lceil \frac{2\sigma^2}{\epsilon^2} \rceil, n \big\}}$ and any $b' \leq \sqrt{b}$, then the number of iterations can be bounded by 
		\vspace{-1mm}
		$$ T= \frac{4 L \fgap }{\epsilon^2} \left( 1+ \sqrt{\frac{1-p}{pb'}}\right) + \frac{1}{p} \leq \frac{8L\fgap \sqrt{b}}{\epsilon^2 b'} + \frac{b+b'}{b'},
		$$
		and the number of stochastic gradient computations (i.e., gradient complexity) is 
		$$ \#\mathrm{grad} =b+ T\left(pb+(1-p)b'\right) \leq  3b + \frac{16L\fgap \sqrt{b}}{\epsilon^2} = \conline{O\bigg(b+\frac{\sqrt{b}}{\epsilon^2}\bigg)}.
		$$
	\end{enumerate}
\end{theorem}
}	\vspace{-5mm}

\section{Simple Convergence Analysis for Theorem \ref{thm:main}}
\label{sec:proof}
First, we use the following key Lemma \ref{lem:relation} which describes a useful relation between the function values after and before a gradient descent step, i.e.,  between $f(\xtn)$ and $f(\xt)$ with $\xtn := \xt - \eta \gt$ for any gradient estimator $g^t \in \R^d$ and stepsize $\eta >0$.
\begin{lemma}[Lemma 2 of \citealp{li2021page}]
	\label{lem:relation}
	Suppose that function $f$ is $L$-smooth and let $\xtn := \xt - \eta \gt$. Then for any $g^t\in \R^d$ and $\eta>0$, we have
	\begin{align}\label{eq:relation}
	f(\xtn) \leq f(\xt) - \frac{\eta}{2} \ns{\nabla f(\xt)} 
	- \Big(\frac{1}{2\eta} - \frac{L}{2}\Big) \ns{\xtn -\xt}
	+ \frac{\eta}{2}\ns{\gt - \nabla f(\xt)}.
	\end{align}
\end{lemma}

\vspace{1mm}
\topic{\cfinite{Finite-sum} case} Then we use the following Lemma \ref{lem:var-finite} to deal with the last variance term of \eqref{eq:relation} for the {finite-sum} case \eqref{prob:finite}.
\begin{lemma}[Lemma 3 of \citealp{li2021page}]
	\label{lem:var-finite}
	Suppose that Assumption \ref{asp:smooth} holds. If the gradient estimator $\gtn$ is defined in Line \ref{line:grad} of Algorithm \ref{alg:page}, then we have
	\begin{align}
	\E\big[\ns{\gtn-\nabla f(\xtn)}\big] \leq  (1-p_t) \ns{g^{t} - \nabla f(\xt)} + \frac{(1-p_t) L^2}{b'}\ns{\xtn - \xt}. \label{eq:var-finite} 
	\end{align}
\end{lemma}

Now, we are ready to prove Theorem \ref{thm:main} by combining Lemmas \ref{lem:relation} and \ref{lem:var-finite}. 
We add \eqref{eq:relation} with $\frac{\eta}{2p}$ $\times$ \eqref{eq:var-finite} (here we simply let $p_t \equiv p$), and take expectation to get
\begin{align}
&\EB{ f(\xtn) -f^* + \frac{\eta}{2p}\ns{\gtn-\nabla f(\xtn)}}  \notag \\
&\leq \EB{ f(\xt) -f^* - \frac{\eta}{2} \ns{\nabla f(\xt)} 
	- \Big(\frac{1}{2\eta} - \frac{L}{2}\Big) \ns{\xtn -\xt}
	+ \frac{\eta}{2}\ns{\gt - \nabla f(\xt)}} \newl
+ \frac{\eta}{2p} \EB{(1-p) \ns{g^{t} - \nabla f(\xt)} + \frac{(1-p) L^2}{b'}\ns{\xtn - \xt}} \notag\\
&= \EB{f(\xt) -f^* + \frac{\eta}{2p}\ns{\gt-\nabla f(\xt)}  - \frac{\eta}{2} \ns{\nabla f(\xt)} 
	-\Big(\frac{1}{2\eta} - \frac{L}{2} -\frac{(1-p)\eta L^2}{2pb'}\Big) \ns{\xtn -\xt}} \notag\\
&\leq \EB{f(\xt) -f^* + \frac{\eta}{2p}\ns{\gt-\nabla f(\xt)}  - \frac{\eta}{2} \ns{\nabla f(\xt)}}, \label{eq:use-eta}  
\end{align}
where $f^*:=\min_{x\in \R^d}f(x)$ and the last inequality \eqref{eq:use-eta} holds due to $\frac{1}{2\eta} - \frac{L}{2} -\frac{(1-p)\eta L^2}{2pb'} \geq 0$ by choosing stepsize 
\begin{align}\label{eq:eta}
\eta \leq \frac{1}{L\Big(1+\sqrt{\frac{1-p}{pb'}}\Big)}.
\end{align} 
If we define $\Phi_t := f(\xt) -f^* + \frac{\eta}{2p}\ns{\gt-\nabla f(\xt)}$, then 
summing up \eqref{eq:use-eta} from $t=0$ to $T-1$, we get
\begin{align}
\E[\Phi_T] \leq \E[\Phi_0] - \frac{\eta}{2} \sum_{t=0}^{T-1}\E[\ns{\nabla f(\xt)}]. \label{eq:phit}
\end{align} 
Thus according to the output of \page, i.e.,  $\hx_T$ is randomly chosen from $\{\xt\}_{t\in[T]}$ and $\Phi_0:=f(x^0) -f^* + \frac{\eta}{2p}\ns{g^0-\nabla f(x^0)}=f(x^0) -f^* \overset{\text{def}}{=} \fgap$, we have
\begin{equation}\label{eq:final}   
\E[\ns{\nabla f(\hx_T)}] \leq \frac{2\fgap}{\eta T} = \epsilon^2,
\end{equation}
where the last equality holds if the number of iterations is 
\begin{align}
T= \frac{2\fgap}{\epsilon^2 \eta}  
\overset{\eqref{eq:eta}}{=}  \frac{2\fgap L}{\epsilon^2}  \bigg(1+\sqrt{\frac{1-p}{pb'}}\bigg).
\end{align}
Note that $
\E[\n{\nabla f(\hx_T)}] \leq \sqrt{\E[\ns{\nabla f(\hx_T)}]} \leq \epsilon
$ from \eqref{eq:final} via Jensen's inequality.
\qedb

\newpage
\head{\conline{Online} case} Similarly, for the {online} case \eqref{prob:online}, we use the following Lemma \ref{lem:var-online} instead of Lemma \ref{lem:var-finite} to deal with the last variance term of \eqref{eq:relation}. 
Note that we refer the online problem \eqref{prob:online} as the finite-sum problem \eqref{prob:finite} with large or infinite $n$.
The additional bounded variance Assumption \ref{asp:bv} usually is needed for this online case since the full gradient ($\nabla f(x)=\frac{1}{n}\sum_{i=1}^n{\nabla f_i(x)}$) may not be available.

\begin{lemma}[Lemma 4 of \citealp{li2021page}]
	\label{lem:var-online}
	Suppose that Assumptions \ref{asp:smooth} and \ref{asp:bv} hold. If the gradient estimator $\gtn$ is defined in Line \ref{line:grad} of Algorithm \ref{alg:page}, then we have
	\begin{align}
	\E[\ns{\gtn-\nabla f(\xtn)}] \leq  (1-p_t) \ns{g^{t} - \nabla f(\xt)} + \frac{(1-p_t) L^2}{b'}\ns{\xtn - \xt} + {\bf 1}_{\{b<n\}} \frac{p_t\sigma^2}{b}, \label{eq:var-online} 
	\end{align}
	where ${\bf 1_{\{\cdot\}}}$ denotes the indicator function.
\end{lemma}

Then the remaining proof is similar to previous \cfinite{finite-sum} case. 
Similar to \eqref{eq:use-eta}, here we add \eqref{eq:relation} with $\frac{\eta}{2p}$ $\times$ \eqref{eq:var-online} to get
\begin{align}
&\EB{ f(\xtn) -f^* + \frac{\eta}{2p}\ns{\gtn-\nabla f(\xtn)}}  \notag \\
&\leq \EB{f(\xt) -f^* + \frac{\eta}{2p}\ns{\gt-\nabla f(\xt)}  - \frac{\eta}{2} \ns{\nabla f(\xt)}  
	+ {\bf 1}_{\{b<n\}} \frac{\eta\sigma^2}{2b}}, \label{eq:use-eta-online}  
\end{align}
where the stepsize $\eta$ is chosen the same as in \eqref{eq:eta}. 
Using the same definition $\Phi_t := f(\xt) -f^* + \frac{\eta}{2p}\ns{\gt-\nabla f(\xt)}$ and 
summing up \eqref{eq:use-eta-online} from $t=0$ to $T-1$, we get
\begin{align}\label{eq:last1}
\E[\Phi_T] \leq \E[\Phi_0] - \frac{\eta}{2} \sum_{t=0}^{T-1}\E[\ns{\nabla f(\xt)}]+{\bf 1}_{\{b<n\}} \frac{\eta T \sigma^2}{2b}.
\end{align} 
For the term $\E[\Phi_0]$, we have 
\begin{align}
\E[\Phi_0] :=\E\Big[f(x^0) -f^* + \frac{\eta}{2p}\ns{g^0-\nabla f(x^0)}\Big]
&=\E\bigg[ f(x^0) -f^*  + \frac{\eta}{2p}\Big\|\frac{1}{b} \sum_{i\in I} \nabla f_i(x^0) -\nabla f(x^0)\Big\|^2\bigg] \label{eq:use-g0} \\
&\leq f(x^0) -f^* + {\bf 1}_{\{b<n\}} \frac{\eta \sigma^2}{2pb}, \label{eq:use-b0} 
\end{align}
where \eqref{eq:use-g0} follows from the definition of $g^0$ (see Line \ref{line:sgd} of Algorithm \ref{alg:page}), and 
\eqref{eq:use-b0} is due to Assumption~\ref{asp:bv}, i.e., \eqref{eq:bv}.
Plugging \eqref{eq:use-b0} into \eqref{eq:last1} (noting that $\fgap:= f(x^0) -f^*$) and according to the output of \page, i.e.,  $\hx_T$ is randomly chosen from $\{\xt\}_{t\in[T]}$, we have 
\begin{align}
\E[\ns{\nabla f(\hx_T)}] &\leq \frac{2\fgap}{\eta T} +  {\bf 1}_{\{b<n\}}  \frac{\sigma^2}{pbT}  + {\bf 1}_{\{b<n\}} \frac{\sigma^2}{b} \notag\\
&\leq \Big(\frac{2\fgap}{\eta} +  \frac{\epsilon^2}{2p}\Big) \frac{1}{T}  + \frac{\epsilon^2}{2} \label{eq:chooseb}\\
& =\epsilon^2,\label{eq:last-online}
\end{align}
where \eqref{eq:chooseb} follows from the parameter setting of minibatch size $b=\min\big\{\lceil \frac{2\sigma^2}{\epsilon^2} \rceil, n\big\}$, and the last equality \eqref{eq:last-online} holds if the number of iterations 
\begin{align}
T= \frac{4\fgap}{\epsilon^2 \eta} + \frac{1}{p}   
\overset{\eqref{eq:eta}}{=}  \frac{4\fgap L}{\epsilon^2}  \left(1+\sqrt{\frac{1-p}{pb'}}\right) + \frac{1}{p}.
\end{align}
\qedb


\bibliographystyle{plainnat}
\bibliography{note_page}

\end{document}